\documentclass[10pt]{article}

\usepackage{latexsym}
\usepackage{amssymb}
\usepackage{amsmath}
\usepackage{amsthm}
\usepackage{mathrsfs}
\input{epsf}

\def\AFOUR{%
\setlength{\textheight}{9.0in}%
\setlength{\textwidth}{5.75in}%
\setlength{\topmargin}{-0.375in}%
\hoffset=-.5in%
\renewcommand{\baselinestretch}{1.17}%
\setlength{\parskip}{6pt plus 2pt}%
}

\AFOUR                                           

\newtheorem{theorem}{Theorem}[section]
\newtheorem{lemma}[theorem]{Lemma}
\theoremstyle{definition}
\newtheorem{definition}[theorem]{Definition}

\theoremstyle{remark}
\newtheorem{remark}[theorem]{Remark}
\newtheorem{proposition}[theorem]{Proposition}

\numberwithin{equation}{section}

\newcommand{\nc}{\newcommand}
\newcommand{\rnc}{\renewcommand}

\nc{\bea}{\begin{eqnarray}}
\nc{\eea}{\end{eqnarray}}
\nc{\be}{\bea}
\nc{\ee}{\eea}

\nc{\ex}[1]{\mbox{e}^{\,\textstyle#1}}

\nc{\eps}{\epsilon}
\nc{\tr}{\mathop{\mbox{tr}}\nolimits}
\nc{\ad}{\mathop{\mbox{ad}}\nolimits}
\nc{\Tr}{\mathop{\mbox{Tr}}\nolimits}
\nc{\Det}{\mathop{\mbox{Det}}\nolimits}
\nc{\rk}{\mathop{\mbox{rk}}\nolimits}

\rnc{\lg}{\mathbf{g}}
\nc{\lt}{\mathfrak{t}}

\begin{document}
\global\parskip=4pt
\begin{center}

{\Large\bf Intersection
  Pairings on Spaces of Connections 
 and Chern-Simons Theory on Seifert Manifolds}
\vskip 0.4in
{\bf George Thompson}\\

ICTP, Strada Costiera 11, 34100 Trieste, Italy.
\date{}
\end{center}

\section{Introduction}
Let $M$ be a compact closed 3-manifold and $G \rightarrow P
\rightarrow M$ a trivial principal $G$ bundle over $M$. We denote by
$\mathcal{A}_{3}$ the space of connections on $P$ and $\mathcal{G}_{3}$ the
gauge group. Let $Z_{k,G}[M]$ denote the Chern-Simons path integral at level $k$
and for group $G$ (taken as above),
\be
Z_{k,G}[M]= \frac{1}{\mathrm{Vol}(\mathcal{G}_{3})} \int_{\mathcal{A}_{3}}\,
  \exp{\left(  I(\mathscr{A})\right)} \nonumber
\ee
where the Chern-Simons action, for $\mathscr{A} \in \mathcal{A}_{3}$, is
\be
I(\mathscr{A})= i\frac{k}{4\pi}\int_{M} \Tr{\left(\mathscr{A} \wedge
    d\mathscr{A}  +
    \frac{2}{3}\mathscr{A} \wedge \mathscr{A}  \wedge \mathscr{A}
  \right)} \nonumber 
\ee
and $\Tr$ is normalized so that under large gauge transformations
(recall that $\pi_{0}(\mathcal{G}_{3}) = \mathbb{Z}$)
$I(\mathscr{A}^{\mathscr{G}})=I(\mathscr{A}) + 2\pi i n$, $\mathscr{G}
\in \mathcal{G}_{3}$ and $n \in \mathbb{Z}$ so that the exponential is
invariant.  

One can also consider the inclusion of `Wilson lines' in the path
integral. A Wilson line is a
combination of a knot $K$ and an irreducible representation $R$ of the group
$G$, and is defined to be the holonomy
\be
W_{R}(K) = \Tr_{R}{P\exp{\oint_{K}\mathscr{A}}} \nonumber
\ee
The relevant path integral is
\be
Z_{k,G}[M,(K_{i},R_{i})] = \frac{1}{\mathrm{Vol}(\mathcal{G}_{3})}
\int_{\mathcal{A}_{3}}  \exp{\left(I(\mathscr{A}) \right)}\, 
\prod_{i=1} W_{R_{i}}(K_{i}) \nonumber
\ee
Let $\Sigma$ be a smooth genus $g$ curve and $\omega$ a unit
volume K\"{a}hler form on $\Sigma$. Let $\mathfrak{M}$ be the
moduli space of flat $G$ connections on $\Sigma$. The moduli
space has a natural symplectic, infact K\"{a}hler, structure and we
let $\mathfrak{L}
\rightarrow \mathfrak{M}$ be the fundamental line bundle 
whose first Chern class agrees with the natural symplectic
form on $\mathfrak{M}$ (for $G=SU(n)$ this will be the
determinant line bundle). E. Witten \cite{Witten-CS} has
shown that the quantum Hilbert space of 
states of Chern-Simons theory on $\Sigma$ is
$\mathrm{H}^{0}(\mathfrak{M},
\mathfrak{L}^{k})$. Quite generally the 
Chern-Simons invariant, $Z_{k,G}[M]$, for a 3-manifold $M$, by Heegard
splitting, will be the inner product of two vectors in the Hilbert
space. However, the dimension of the Hilbert space is the Chern-Simons
invariant of the 3-manifold 
$\Sigma \times S^{1}$. 

If one includes Wilson lines, then the
Hilbert space is $\mathrm{H}^{0}(\mathfrak{M}, \mathfrak{L}^{k})$ where now
$\mathfrak{M}$ is the moduli space of parabolic bundles on
$\Sigma$, which is still naturally K\"{a}hler, and $\mathfrak{L}$ is
the associated fundamental line bundle.

The
Hirzebruch-Riemann-Roch theorem tells us that
\be
\sum_{q=0} (-1)^{q}\mathrm{dim}\,\mathrm{H}^{q}(\mathfrak{M},
\mathfrak{L}^{k}) 
= \int_{\mathfrak{M}} \mathrm{Todd}(\mathfrak{M}) \wedge
\mathrm{Ch}(\mathfrak{L}^{k} ) \nonumber
\ee
If the canonical bundle of $\mathfrak{M}$ is negative, as will be
the case in the examples we discuss, then the higher cohomology groups
are trivial by Kodaira vanishing and we have
\be
\mathrm{dim}\,\mathrm{H}^{0}(\mathfrak{M}, \mathfrak{L}^{k})
= \int_{\mathfrak{M}} \mathrm{Todd}(\mathfrak{M}) \wedge
\mathrm{Ch}(\mathfrak{L}^{k} ) \nonumber 
\ee
E. Verlinde \cite{Verlinde} provides us with a
concrete formula for the dimension of $\mathrm{H}^{0}(\mathfrak{M},
\mathfrak{L}^{k})$. However, as we have seen, Chern-Simons theory
provides us with another formulation of 
E. Verlinde's dimension count, namely
\be
Z_{k,G}[\Sigma \times S^{1},(K_{i},R_{i})] = \int_{\mathfrak{M}}
\mathrm{Todd}(\mathfrak{M}) \wedge 
\mathrm{Ch}(\mathfrak{L}^{k} ) \nonumber
\ee
The moduli spaces, which thus far have been generically denoted by
$\mathfrak{M}$, are singular at the points where there are reducible
connections. Nevertheless, suitably interpreted, the index theorem
still yields a topological expression for the invariants $Z_{k,G}[\Sigma \times
S^{1},(K_{i},R_{i})]$. This raises a

\vspace{0.3cm}
{\bf Question:}
Are there other
3-manifolds whose Chern-Simons invariants, or parts of them, can be expressed as
intersection pairings on an appropriate moduli space $\mathfrak{M}$?
\\

The question has been partially answered in the affirmative by Beasley and
Witten \cite{Beasley-Witten}. From now on set  $M \equiv M_{(g,p)}$ to
denote the Seifert manifold
that is presented as a degree $-p$, $U(1)$ bundle over a Riemann
surface $\Sigma$ of genus $g$.
Using non-Abelian localization for
the Seifert 3-manifolds $M$ Beasley and Witten are
able to show that (equation (5.176) in their paper with $n=-p$ due to a
different choice of orientation and noting that there is a slightly
different normalization of $\Theta$),
\begin{proposition}\label{BW}{{\rm (Beasley-Witten) Let $M$ be as above. The
      portion of the Chern-Simons invariant which is localized on the
      smooth part, $\mathfrak{M}$, of the moduli space of Yang-Mills
      connections is 
\be
\left. Z_{k, \, G}[M]\right|_{\mathfrak{M}} 
 = 
\frac{1}{|\Gamma|}\exp{\left(i\frac{\pi}{2} \eta_{0}\right)}\,
\int_{\mathfrak{M}} \mathrm{Todd}(\mathfrak{M})\wedge
\mathrm{Ch}(\mathfrak{L}^{k} ) \wedge \exp{\left(
    - i \frac{p}{2\pi}(k+c_{\lg})
    \Theta(\mathfrak{M}) \right) } \nonumber
\ee
In this formula, $\Theta(\mathfrak{M})$ is a certain degree 4
cohomolgy class on $\mathfrak{M}$, $c_{\lg}$ and $\Gamma$ are the dual Coxeter
number and centre of $G$ respectively while $\eta_{0}$ is the framing
of $M$. Notice that this formula is quite analogous to the
Hirzebruch-Riemann-Roch formula and reduces
to it when $p=0$.
}}
\end{proposition}
Let $\mathcal{A}$ be the space of connections on the trivial $G$ bundle on
$\Sigma$ and $\mathcal{G}$ the associated group of gauge
transformations. I will show that
\begin{proposition}\label{onS}{{\rm  The Chern-Simons path integral on
$M$, $Z_{k, \, G}[M]$, is equal to a path integral on the space of
connections over $\Sigma$, namely
\be
Z_{k, \, G}[M] = \exp{\left(i\frac{\pi}{2}
      \eta_{0}\right)}\frac{1}{\mathrm{Vol}(\mathcal{G})} 
\int_{\mathcal{A}} \mathrm{Todd}(\mathcal{A})\wedge
\mathrm{Ch}(\mathfrak{L}^{k} ) \wedge \exp{\left(
    - i \frac{p}{2\pi}(k+c_{\lg})
    \Theta(\mathcal{A}) \right) } \nonumber
\ee
in this formula all the classes have the same form as those in
Proposition \ref{BW}.
}}
\end{proposition}
As there is a striking resemblance between the formulae presented in
Propositions \ref{BW} and \ref{onS} it is worthwhile, at this point, to
make some remarks.  
\begin{remark}{{\rm The formula in Proposition \ref{onS} is an exact
      expression for $Z_{k, \, G}[M]$ unlike that in Proposition
      \ref{BW} which is a part of the answer. 
}}
\end{remark} 
\begin{remark}{{\rm It is quite straightforward to evaluate the path
      integral in Proposition \ref{onS}. This has been done for
      general $\Sigma$ in
      \cite{Blau-Thompson-CS} by Abelianization and for $\Sigma=S^{2}$
      in \cite{Beasley-Witten} by non-Abelian localization. The
      Reshetikhin-Turaev-Witten 
      invariants for Seifert 3-manifolds may be 
obtained by surgery prescriptions
as in \cite{Lawrence-Rozansky} and \cite{Marinio} and agree with the
path integral results. However, I leave the path integral
`un-integrated' as it brings the geometry to the fore.
}}
\end{remark}
\begin{remark}{{\rm I am not at all implying that Proposition \ref{BW}
      follows easily from Proposition \ref{onS}, though one would reasonably
      expect that it does. Presumably an application of non-Abelian
      localization to the path integral appearing in Proposition
      \ref{onS} is what is required.
}}
\end{remark}
\begin{remark}{{\rm Proposition \ref{onS} has already been
      established, somewhat indirectly in \cite{AOSV} and more
      directly in \cite{Blau-Thompson-CS}, though the language is somewhat
      different and the classes were not identified in either of these
      works and so the geometric significance of the right hand side
      of Proposition \ref{onS} was not appreciated.
}}
\end{remark}
I will prove a rather more general result involving knots in the fibre
direction which are located at points $x_{i}\in \Sigma$ on the base of
the fibration and which run along the fibre. The point of view adopted here
is that to the parabolic points $x_{i}$ one `attaches' a co-adjoint
orbit $M_{R_{i}}$ defined by the representation $R_{i}$.
\begin{proposition}\label{onSKnot}{{\rm The Chern-Simons path integral on
$M$, with knots in the fibre direction at the points $x_{i} \in
\Sigma$ and associated representations $R_{i}$, $Z_{k, \, G}[M,
(x_{i},R_{i})]$, is given by
\bea
Z_{k,G}[M,(x_{i},R_{i})] 
 & = & 
\exp{\left(i\frac{\pi}{2}\eta_{0}\right)}\frac{1}{\mathrm{Vol}(\mathcal{G})}\,
\int_{\mathcal{A}\times   \prod_{i}   M_{R_{i}}} \widehat{A}(\mathcal{A}\times
\prod_{i}M_{R_{i}}) \, \wedge  \nonumber\\
& &\;\;\;\;\; .\, \exp{\left( (k+c_{\lg})\Omega(\mathcal{A})
    +\sum_{i}\omega(M_{R_{i}}) 
    -i \frac{p}{2\pi}(k+c_{\lg}) \Theta(\mathcal{A}) \right)
}\nonumber
\eea
where $\widehat{A}$ is the $A$ hat genus.
}}
\end{proposition}

{\bf Acknowledgements} This paper is an outcome of a line of research
on abelianization in low dimensions that Matthias Blau and I have been
pursuing over many years. It is a pleasure to thank him for this very
enjoyable collaboration. It is also a pleasure to thank
M.S. Narasimhan and T. Ramadas for helping me to formulate the results
presented here. Finally a thank you to the organizers
of the Bonn workshop ``Chern-Simons Gauge Theory: 20 Years After'' and
to Chris Beasely for discussions on the relationship between his work
and mine.

\section{Background and Strategy}

While it is clear from the description of the Hilbert space of
Chern-Simons theory that the moduli spaces are built in, what is not
at all obvious is that Chern-Simons should, in general, have knowledge of the
cohomology ring of the moduli space. This paragraph is intended to
motivate such a connection. Recall that Witten \cite{Witten-YM2Rev}
had established that the topological field theory analogue of
Donaldson theory on a curve can be mapped to Yang-Mills theory on the
curve. Since the topological field theory is designed to probe the
cohomology ring of the moduli space we learn that Yang-Mills theory
will do just that. The
action of Yang-Mills theory is 
\be
S(F_{A},\psi,\phi) = \frac{1}{4\pi^{2}}\int_{\Sigma} \Tr{\left(i\phi
  F_{A} + \frac{1}{2} \psi \wedge \psi \right)} +
\frac{\eps}{8\pi^{2}} \int_{\Sigma} \omega
\,\Tr{\phi^{2}} \label{YM2Action} 
\ee
where $A$ is a connection on $P$, $\phi \in \Gamma(\Sigma, \ad \, P)$
and $\psi$ is interpreted as a one form on the space
$\mathcal{A}/\mathcal{G}$ (in terms of the universal
bundle construction described in the 
next section the elements of the action are all cohomology classes on
$\mathcal{A}/\mathcal{G}$).

In \cite{Beasley-Witten} Beasley and Witten recall that every
3-manifold has a contact structure $\kappa$. Here we use the $U(1)$
bundle structure and the associated nowhere vanishing vector
field. (This has the advantage that it corresponds to the obvious
structure on $\Sigma \times S^{1}$.) One may, therefore,
decompose connections as
\be
\mathscr{A}=  A + \kappa\,\frac{1}{2\pi} \phi , \;\;\; \iota_{\kappa}A=0 \nonumber
\ee
The Chern-Simons action, on choice of the contact structure, is
\be
I(A,\phi) = i\frac{k}{4\pi^{2}}\int_{M}\left(\pi \kappa \Tr{A \iota_{\kappa} dA}
  +\kappa \Tr{\phi F_{A}} + 
  \kappa d\kappa \, \frac{1}{4\pi}\Tr{\phi^{2}} \right) \label{CS-Action}
\ee
and this has some resemblance to (\ref{YM2Action}). The first
term in (\ref{CS-Action}) is the `time' derivative (in the direction
of the vector field dual to $\kappa$). There are some differences and
perhaps the most glaring is that the term proportional to $\psi \wedge
\psi$ is absent. However, a variation on the theme is now
available. The structure group of a
general 3-manifold is $SO(3)$, however, on fixing the bundle structure
the compatible structure group reduces to $SO(2)$ (the structure group
of $\Sigma$). There is a minimal, so called $N=1$, supersymmetric version
of Chern-Simons theory which is not obviously topological. However,
when the structure group is reduced to $SO(2)$ one may `twist' the
$N=1$ supersymmetric Chern-Simons theory to obtain a manifestly topological
theory (up to choice of contact structure). The twisted version has 
action (\ref{CS-Action}) augmented with
\be
\frac{k}{8\pi^{2}}\int_{M} \kappa\, \Tr{ \psi \wedge \psi}, \;\;\;
\iota_{\kappa}\psi =0 \nonumber
\ee
The $N=1$ action becomes
\be
I(A, \phi, \psi) = \frac{k}{4\pi^{2}}\int_{M}\left(i\pi \kappa
  \Tr{A\wedge \iota_{\kappa}dA } + \kappa \Tr{\left(i \phi F_{A} +
      \frac{1}{2} \psi \wedge \psi \right)} + \frac{i}{4\pi} \kappa
  d\kappa \Tr{\phi^{2}}\right) \nonumber
\ee
so that now the resemblance of the two theories is rather
remarkable. The supersymmetry transformations are closely linked to 
those in \cite{Beasley-Witten}, equation (3.48),
\be
QA= i\psi, \;\;\; Q\psi = - d_{A} \phi - 2 \pi \iota_{\kappa}dA \label{susy}
\ee
The fact that the fields $\psi$ are free means that partition function
and knot invariants of the $N=1$ and the usual Chern-Simons theory agree.

One possible novelty that there is a
new set of observables in the theory namely those that are 3-forms on
the moduli space
\be
\int_{M} J_{C} \wedge \Tr{\phi.\psi} \nonumber
\ee
with $J_{C}$ a de Rham current with delta function support on the
1-dimensional cycle $C \subset M$ is, unfortunately, spoilt by the
second part of the transformation of $\psi$ in (\ref{susy}).

In any case it would appear that Chern-Simons theory is some
supersymmetric quantum mechanics on the moduli spaces of
interest. Since supersymmetric quantum mechanics is related to index
theory one may reasonably expect a positive response to the question.

The strategy we will follow here is to
\begin{enumerate}
\item Identify the objects that appear in the path integral in terms
  of cohomology classes on the space of connections modulo the gauge group. This
  is done in Section \ref{UniversalBundle}.
\item Next we relate the classes that one obtains from the universal bundle
  construction to classes associated with the tangent bundle of
  $\mathcal{A}/\mathcal{G}$. These are denoted by
  $\mathrm{Todd}(\mathcal{A})$ and 
  $\widehat{A}(\mathcal{A})$ in Section \ref{ToddAhat} for reasons
  that will become apparent as we go along.
\item One can now express Witten's localisation of the Yang-Mills path
  integral on $\Sigma$ (Proposition
  \ref{wittenloc} below) in terms of classes that come from the
  universal bundle on $\mathcal{A}/\mathcal{G}$ on the one hand and
  the same classes restricted to $\mathfrak{M}$ on the other.
\item Having established the relationship with the cohomology ring in
  $\mathcal{A}/\mathcal{G}$ our final task is to show that the
  integral over $\mathcal{A}_{3}/\mathcal{G}_{3}$ reduces to that of certain
  classes on $\mathcal{A}/\mathcal{G}$ on integrating out sections
  which are not $U(1)$ invariant.
\end{enumerate}

\section{ A Universal Bundle}\label{UniversalBundle}
Let $P$ be a principal  $G$ bundle over some smooth manifold $X$,
$\mathcal{A}$ the space of 
connections on $P$ and $\mathcal{G}$ the group of gauge
transformations (bundle automorphisms). We have the following
universal bundle construction due to 
Atiyah and Singer \cite{Atiyah-Singer}. There is an action of
$\mathcal{G}$ on $P$ so that 
we may form the space (this is not smooth unless one makes further
assumptions and as it stands it is a stack)
\be
\mathcal{Q} = \mathcal{A} \times_{\mathcal{G}} P  \nonumber
\ee 
Now $G$ operates on $\mathcal{Q}$ and
infact $\mathcal{Q}$ is itself the total space of a principle bundle
\be
\mathcal{Q}\rightarrow \mathcal{Q}/G = \mathcal{A}/\mathcal{G}
\times X \nonumber
\ee
There is a natural connection on $\mathcal{Q}$ and from it we can
define a curvature 2-form and then Chern classes, via Chern-Weil
theory, for an associated rank $n$ vector bundle $\mathcal{E} = \mathcal{Q}
\times_{G} \mathbb{C}^{n}\rightarrow \mathcal{A}/\mathcal{G} \times
X$. Finally we restrict to some (hopefully smooth) moduli space
$\mathfrak{M} \subset \mathcal{A}/\mathcal{G}$, for which the above
construction makes sense. 
\begin{remark}{{\rm A detailed description of this bundle and
its relationship to topological field theories can be found in
\cite{Birmingham-Blau-Rokowski-Thompson} chapter 5, especially
sections 5.1 and 5.3.
}}
\end{remark}
\begin{remark}{{\rm The path integral we need to deal with is the one over all
      of $\mathcal{A}/\mathcal{G}$ and not some smooth finite
      dimensional subspace $\mathfrak{M}$. Consequently, we will be
      integrating over the stack.
}}
\end{remark}
From now on we take $X$ to be $\Sigma$. Decompose the curvature
2-form on $\mathcal{E}$ into its Kunneth components as 
\be
1\otimes F_{A} + \Psi +\Phi \otimes 1 \in
\mathrm{H}^{2}(\mathcal{A}
\times \Sigma, \mathrm{Lie}\, G) \nonumber
\ee
(If $\mathfrak{M}$ is simply connected we have that $\Tr{\Psi}$,
restricted to $\mathfrak{M}$, is
cohomologously trivial in which case $
c_{1}(\mathcal{E}) = 1 \otimes c_{1}(E) +\frac{i}{2\pi}
\Tr{\Phi}\otimes  1$
where $\Tr \Phi \in \mathrm{H}^{2}(\mathfrak{M})$.) If we fix on
$G=SU(n)$ then $c_{1}(\mathcal{E})$ vanishes and
the second Chern class decomposes as
\bea
 c_{2}(\mathcal{E}) & = &  
\frac{1}{4\pi^{2}}\Tr{\left( \Phi \otimes 
    F_{A} +\frac{1}{2} \Psi \wedge \Psi \right) } + \frac{1}{4\pi^{2}}\Tr{
  \Phi \wedge \Psi}+
\frac{1}{8\pi^{2}}\Tr{\Phi^{2}} \otimes 1 \nonumber\\
&= & \Omega(\mathcal{E}) + \gamma(\mathcal{E}) +
\Theta(\mathcal{E}) \label{c2kunneth} 
\eea
We have a differential $Q$ on the space $\mathcal{A} \times \Sigma$
which satisfies
\be
QA = \Psi, \;\;\; Q\Psi = d_{A}\Phi, \;\;\; Q\Phi =0 , \;\;\; Q^{2} =
\mathcal{L}_{\Phi} \nonumber
\ee
so that $\Psi$ is a (basis) one form on $\mathcal{A}$ and
\be
(Q-d_{A})(1\otimes F_{A} + \Psi +\Phi \otimes 1) = 0 .\nonumber
\ee
Hence associated to $Q$ we have equivariant cohomology on
$\mathcal{A}/\mathcal{G} \times \Sigma$ and the Chern classes
$c_{n}(\mathcal{E})$ are $Q$ closed.

Note the second Chern class appears as the
action of Yang-Mills theory (\ref{YM2Action}) if we make the
identifications $\Phi = i \phi$ and $\Psi = \psi$, and we consider
these as forms on $\mathcal{A}/\mathcal{G} \times \Sigma$
\be
S(F_{A},\psi,\phi) = \pi_{*} \left(\Omega(\mathcal{E}) - \eps
\Theta(\mathcal{E})\otimes \omega \right)\simeq
\Omega(\mathcal{A})-\eps \Theta(\mathcal{A}) \label{c2isym}
\ee
where $\pi:
\mathcal{A}/\mathcal{G} \times \Sigma
\rightarrow \mathcal{A}/\mathcal{G}$ is projection
onto the first factor.
\begin{remark}{{\rm The identification (\ref{c2isym}) shows us that the
      Yang-Mills action is a $Q$ closed form (of mixed degree). It is
      also quite clearly not $Q$ exact.
}}
\end{remark}

In general for a vector
bundle $V$ $c_{2}(\mathrm{End}\, V) = 2rc_{2}(V) - (r-1)c_{1}(V)^{2}$ and
consequently $c_{2}(\mathrm{End}\, \mathcal{E})
=2rc_{2}(\mathcal{E})$. The main interest here will be on the trace 
free part $\mathrm{End}_{0}\, \mathcal{E}$ (which in any case has the
same second Chern class as $\mathrm{End}\, \mathcal{E}$). Note that 
the classes on $\mathcal{E}$ are in the `fundamental' representation
while they are taken to be in the `adjoint' representation on
$\mathrm{End}_{0}\,\mathcal{E}$.

In the rank 2 case, Newstead \cite{Newstead},
writes the second Chern Class of his universal bundle as
\be
c_{2}(\mathrm{End}\, U)=2\alpha \otimes \omega + 4 \gamma -
  \beta\otimes 1 
\ee
Thus to make contact with that work restrict to
$\mathfrak{M}\subset \mathcal{A}/\mathcal{G}$ and set,
\be
\alpha = \frac{1}{4\pi^{2}} \int_{\Sigma}\Tr{ \Psi \wedge \Psi}, \;\;\;
\beta = -\frac{1}{2\pi^{2}}\Tr{\Phi^{2}}  , \;\;\;
\gamma = \frac{1}{4\pi^{2}}\Tr{
  \Phi \wedge \Psi} \nonumber
\ee

\section{The Todd and the $\widehat{A}$ Genera}\label{ToddAhat}

In order to express the Todd genus in terms of the classes arising from
the universal bundle we follow the approach of Newstead
\cite{Newstead} for determining the Pontrjagin class in the rank 2
case. Then we use an observation of Thaddeus \cite{Thaddeus} to give
the $\widehat{A}$ class in terms of the Pontrjagin roots and from
there Todd. 

The tangent bundle,
$\mathrm{T}_{\mathfrak{M}}$, of $\mathfrak{M}$ is given by
\be
\mathrm{T}_{\mathfrak{M}} \simeq R^{1} \pi_{*} \mathrm{End}_{0}\,\mathcal{E}
 \nonumber
\ee
where $R^{i}$ denotes the $i$-th direct image sheaf under the map $\pi:
\mathfrak{M} \times \Sigma \rightarrow \mathfrak{M}$ onto
the first factor. The Grothendieck-Riemann-Roch theorem states that
\be
\mathrm{Ch}(\mathrm{T}_{\mathfrak{M}})-
\mathrm{Ch}(\pi_{*}\mathrm{End}_{0}\,  \mathcal{E})  = -
\pi_{*}\left(\mathrm{Ch}(\mathrm{End}_{0}\, 
  \mathcal{E}) (1-(g-1) \omega) \right) \nonumber
\ee
For the spaces that we are interested in the direct image sheaf
$R^{0}\pi_{*}\mathrm{End}_{0}\, \mathcal{E}$ is trivial so
\be
\mathrm{Ch}(\mathrm{T}_{\mathfrak{M}})  = -
\pi_{*}\left(\mathrm{Ch}(\mathrm{End}_{0}\, 
  \mathcal{E}) (1-(g-1) \omega \right)\label{Groth-Riemann-Roch}
\ee

Denote the complexification of the Lie algebra of $G$ by
$\lg_{\mathbb{C}}$, the complexification of the Cartan subalgebra
by $\lt_{\mathbb{C}}$ and the space of roots by $\mathbf{k}$ then
\be
\lg_{\mathbb{C}}=\lt_{\mathbb{C}} \oplus \mathbf{k}\nonumber
\ee
We also set $\mathbf{k}_{+}$ to be the space of positive roots.

\begin{lemma}\label{Pontjagin}{\rm{The Pontrjagin class of the tangent bundle
      $P(T_{\mathfrak{M}})$ 
is given by
\be
P(T_{\mathfrak{M}})= \mathrm{det}_{\mathbf{k}}( 1 + \mathrm{ad}\,
\Phi/2\pi)^{2g-2} = \prod_{\mathbf{k}_{+}} \left(1 -
\left(\frac{\alpha(\Phi)}{2\pi}\right)^{2}\right)^{2g-2}  \nonumber
\ee
}}
\end{lemma}
{\bf Proof:} Quite generally the only classes that
contribute to $\mathrm{Ch}_{2m}(\mathrm{T}_{\mathfrak{M}})$ from
(\ref{Groth-Riemann-Roch}) are powers of $\Phi$, so we have
$
\mathrm{Ch}_{2m}(\mathrm{T}_{\mathfrak{M}})
=(g-1)\Tr_{\mathrm{Ad}}{\left(\exp{\left(i\Phi/2\pi\right)}\right)_{2m}}$,
and consequently $
\mathrm{Ch}(\mathrm{T}_{\mathfrak{M}} \oplus
\mathrm{T}_{\mathfrak{M}}^{*}) = 2(g-1)
\Tr_{\mathrm{Ad}}{\left(\exp{\left(i\Phi/2\pi\right)}\right)}$. This is the same as
having the direct sum of $2g-2$ copies of a vector bundle with curvature 2-form
$\mathrm{ad}\,\Phi$. The Chern class for one copy is $\mathrm{det}_{\lg}( 1
+i \mathrm{ad}\, \Phi/2\pi)= \mathrm{det}_{\mathbf{k}}( 1
+ i\mathrm{ad}\, \Phi/2\pi)$, the direct sum formula for Chern classes
gives us that $c(\mathrm{T}_{\mathfrak{M}}\oplus
\mathrm{T}_{\mathfrak{M}}^{*})= \prod_{\alpha \in \mathbf{k}_{+}}(1+
(\alpha(\Phi)/2\pi)^{2})^{2g-2}$. On the otherhand we have that the
Pontrjagin classes are related to Chern classes by $p_{m}(E)=
(-1)^{m}c_{2m}(E\oplus E^{*})$ and this completes the proof. \qed
\begin{lemma}{\rm{The Todd class of the tangent bundle of
      $\mathfrak{M}$ is
\be
\mathrm{Todd}(\mathfrak{M}) = \exp{\frac{1}{2}c_{1}(T_{\mathfrak{M}})} \,
. \, \left( 
  \frac{\det_{\mathbf{k}}{\sin{(\ad{\Phi}/4\pi)}}}{\det_{\mathbf{k}}{(\ad{
        \Phi}/4\pi})} 
\right)^{1-g} \nonumber
\ee
}}
\end{lemma}
{\bf Proof:} Thaddeus (p147 in \cite{Thaddeus}) notes that on writing,
$P = \prod_{i=1}(1+y_{i})$ where the $y_{i}$ are the Pontrjagin roots,
then $\widehat{A} = \prod_{i=1} (\sqrt{y_{i}}/2/\sinh{\sqrt{y_{i}}/2 })$.
From Lemma \ref{Pontjagin} the roots are $-(\alpha(\Phi)/2\pi)^{2}$ and
they come with a multiplicity of $2g-2$ so that
\be
\widehat{A}(\mathfrak{M}) = \prod_{\alpha \in \mathbf{k}_{+}}
\left(\frac{\sinh{(i\alpha(\Phi)/4\pi)}}{
    i\alpha(\Phi)/4\pi}\right)^{2-2g} \nonumber
\ee
and the standard relation between  $\mathrm{Todd}$ and $\widehat{A}$
completes the proof. \qed

We can also make use of (\ref{Groth-Riemann-Roch}) to determine the
first Chern class of the moduli space. In this case the term
proportional to $(g-1)$ cannot contribute since, as we have seen, it
contributes to even classes. We have
\be
\mathrm{Ch}(\mathrm{End}_{0}\mathcal{E})=
-c_{2}(\mathrm{End}_{0}\mathcal{E}) + \dots = -2r c_{2}(\mathcal{E}) +
\dots \label{myGRR}
\ee
where the ellipses represent higher degree forms.
\begin{theorem}{\rm{(J-M. Drezet and M.S. Narasimhan \cite{Drezet-Narasimhan})
      The
      first Chern class of the tangent bundle of $\mathfrak{M}(r,d)$ the
      moduli space of holomorphic vector bundles of rank $r$ and
      determinant of degree $d$ is
\be
c_{1}(T\mathfrak{M}(r,d))= 2\, \mathrm{g.c.d}(r,d)\,
\Omega(\mathfrak{M}(r,d)) \nonumber
\ee
where $\Omega(\mathfrak{M}(r,d)) $ is the class of the basic line bundle. 
}}
\end{theorem}
In our case $d=0$ so that
$c_{1}(T\mathfrak{M})=2r\Omega(\mathfrak{M})$ and we set, on comparing
with (\ref{myGRR}),
$\Omega(\mathfrak{M})  =
\pi_{*}c_{2}(\mathcal{E})$ that is
\be
\Omega(\mathfrak{M}) = \frac{1}{4\pi^{2}}\int_{\Sigma} \left(i
  \Tr{\phi F_{A} } + \frac{1}{2}\Tr{\Psi \wedge
  \Psi}\right) \label{symp2form} 
\ee
which is the form one would expect on $F_{A}=0$.

\section{Intersection Pairings on Moduli Spaces}

By comparing to the second Chern class of the universal bundle we see
that, quite generally, one should think of $S(F_{A},\psi,\phi)$ as a
form on the space of connections $\mathcal{A}$. With abuse of
notation we denote those classes on $\mathcal{A}$ by the same
symbols as those on $\mathfrak{M}$,
\be
\frac{1}{\mathrm{Vol}(\mathcal{G})} \int_{\mathcal{A} \otimes \Omega^{0}(\Sigma, \lg)}\,
\exp{\left(S(F_{A},\psi,\phi)\right) } \equiv
\int_{\mathcal{A}/\mathcal{G}} \, 
\exp{\left(\Omega(\mathcal{A})- \eps
    \Theta(\mathcal{A})  \right)}
\ee
Expressing the path integral in this way hides certain things, like
the fact that there is a Gaussian integration over the degree four
class or, as the space $\mathcal{A}/\mathcal{G}$ is infinite
dimensional, one cannot expand the exponential out to pick the top
form term. Indeed a detailed analysis of the path integral shows that
some care needs to be exercised in the interpretation of the right
hand side (and on the left). Nevertheless expressing the path integral
in this way is also very suggestive.

E. Witten \cite{{Witten-YM2Rev}} shows that this path integral
essentially devolves to one on the moduli space,
\begin{proposition}\label{wittenloc}{\rm{(E. Witten
      \cite{{Witten-YM2Rev}}) The path 
      integral localises onto the moduli space of flat connections,
\be
\int_{\mathcal{A}/\mathcal{G}} 
\ex{\left(\Omega(\mathcal{A})- \eps
    \Theta(\mathcal{A})  \right)} = \int_{\mathfrak{M}} 
\ex{\left(\Omega(\mathfrak{M})- \eps
    \Theta(\mathfrak{M})  \right)}\, + \; \mathrm{terms}\;
\mathrm{non-analytic} \;\mathrm{in}\; \eps \nonumber
\ee
and the non-analytic terms vanish as $\eps \rightarrow 0$ (provided
$\mathfrak{M}$ is not singular). The
non-analytic terms arise from higher fixed points of the 
action, that is, from non-flat solutions to $d_{A}*F_{A}=0$.
}}
\end{proposition}
The intersection pairings on the moduli space of flat $G$ connections on
$\Sigma$ as presented by Witten \cite{{Witten-YM2Rev}} agree with
those derived by Thaddeus \cite{Thaddeus} in the rank 2 case and
predicts those for higher rank. 
 
\section{Supersymmetric Quantum Mechanics and Passing from
  Chern-Simons to Yang-Mills}
In \cite{Blau-Thompson-CS} one begins with the Chern-Simons theory,
then integrates out modes in the bundle direction, to be left with a
theory on the Riemann surface. This Abelian theory is then solved to
finally provide one with the invariant for the Seifert
manifold. However, it is also pointed out that the same result could
be obtained by considering instead Yang-Mills theory on $\Sigma$ with the
inclusion of some observables one of which (equation (6.8) there)
written in the current notation is 
\be
j_{\lg}(\phi)^{(1-g)} = \widehat{A}(\mathcal{A}) \nonumber
\ee
This observation was critical for the present study.

In this section the approach of \cite{Blau-Thompson-CS} is followed
`half-way' to the point where we have non-Abelian Yang-Mills theory on
$\Sigma$. In this way we are able to obtain the classes on
$\mathcal{A}/\mathcal{G}$ that one must integrate. I will not repeat
the entire calculation but, rather, explain the essential ingredients
especially those which go beyond \cite{Blau-Thompson-CS}. I should
point out that, when needed, I will consider the section $\phi$ to be
momentarily constant on $\Sigma$ and this will simplify the calculation of the
determinants that we will come across presently. After the
determinants are calculated $\phi$ will be allowed to be non-constant
once more. The justification for
this simplification really comes from knowing that had we abelianized
then $\phi$ would be forced to be constant on $\Sigma$ and so consequently we
lose nothing in making this asumption.

\subsection{Fourier Modes and a Gauge Choice}
In order to begin the calculation we note that, as there is a non
degenerate $S^{1}$ action on $M$, we may
decompose all the sections in terms of characters of that action (a
Fourier series). We may therfore write
\bea
\mathcal{A}_{3} & = & \mathcal{A} \oplus \Omega^{0}(\Sigma, \ad{P})
\oplus_{n \neq 0} \Omega^{1}(\Sigma, L^{\otimes -np}\otimes
\ad{P})\oplus_{n \neq 0} 
\Omega^{0}(\Sigma, L^{\otimes -np}\otimes \ad{P})\nonumber \\
\varphi &=& \sum_{n= -\infty}^{\infty} \varphi_{n} , \;\;
\iota_{\kappa}d \varphi_{n} = 
-2\pi i n \varphi_{n}, \;\;\;\varphi_{n} \in \Omega^{*}(\Sigma, L^{\otimes -np}
\otimes \ad{P}) \nonumber
\eea
where $L$ is the line bundle that defines $M$. 

Now note that there is enough gauge
symmetry to impose the condition 
that the section $\phi$ is constant in the fibre direction, that is
$\iota_{\kappa}.d \phi=0$, and we do this. Alternatively put, we make
the identification,
\be
\mathcal{A}_{3}/\mathcal{G}_{3} \simeq\left( \mathcal{A} \oplus
  \Omega^{0}(\Sigma, \ad{P}) 
\oplus_{n \neq 0} \Omega^{1}(\Sigma, L^{\otimes -np}\otimes \ad{P})
\right)/\mathcal{G} 
\nonumber
\ee
There is a caveat
here as the `natural' measures do not coincide since the non-constant
components of the section $\phi$, in
$\Omega^{0}(\Sigma, L^{\otimes -np} \otimes \ad{P})$ are tangent
vectors to the orbit of $\mathcal{G}_{3}$. To correct for this
mismatch one introduces the Faddeev-Popov ghost determinant,
$\Delta_{FP}(\phi)$, which is essentially the ratio of the volume of
the orbit to that of the group. 

In any case the choice of gauge simplifies the path integral immensely. 

\subsection{Integrating over non-trivial characters}
Our aim is to integrate out all those Fourier modes of fields such
that $n\neq0$.
As, by the gauge condition, $\phi$ has no such modes
and the integral over $\psi$ is Gaussian, we concentrate on the
integral of the $A_{n}$ for $n\neq 0$. Note that (with $A_{0}$ denoted
by $A$ again)
\be
I(A,\phi,\psi) = k S(F_{A},\psi, \phi) + \Delta I
\ee
where $S(F_{A},\psi, \phi)$ is the Yang-Mills action with $\eps = ip/2\pi$
and 
\be
\Delta I = \frac{k}{4\pi^{2}}\int_{\Sigma} \sum_{n\neq 0} \Tr{\left(
    A_{-n}\wedge (2\pi n + 
  \ad{\phi})A_{-n} + \psi_{n} \wedge\psi_{-n} \right)}
\ee
Let
\be
\exp{i\Gamma(A, \phi)} = \int \prod_{n\neq
  0}dA_{n}\,d\psi_{n}\;\Delta_{\mathrm{FP}}(\phi) \, 
\exp{i\Delta I}\nonumber
\ee
where the Faddeev-Popov determinant $\Delta_{\mathrm{FP}}$, as I
mentioned above, takes into account the gauge condition on $\phi$.

\begin{definition}{\rm{We will say that two gauge invariant functions
      on $\mathcal{A}$  are
equivalent if their integrals over $\mathcal{A}/\mathcal{G}$ agree and
we will denote that equivalence by $\simeq$.}}
\end{definition}
\begin{proposition}\label{propdet}{\rm{The supersymmetric
      quantum mechanics path integral gives, for $\phi$ valued in the Cartan
      subalgebra, 
\be
\exp{i\Gamma(A,\phi)} \simeq \exp{\left(i\frac{\pi}{2} \eta_{0} \right)}\,
\widehat{A}(i\phi) \wedge \exp{\left(i \frac{c_{\mathbf{g}}}{4\pi^{2}}
    \int_{\Sigma} \Tr{( \phi.F_{A} + \frac{p}{4\pi} \phi^{2}\omega)}
  \right)}  \nonumber
\ee
where $\eta_{0}$ is the framing correction.
}}
\end{proposition} 
{\bf Proof:} As one can see the action is such that the part of the
connection $A^{\mathbf{k}}$, 
only enters in a Gaussian fashion and so 
may easily be integrated out. The integration gives rise to a
determinant which requires regularization (a definition). This
calculation has been performed in \cite{Blau-Thompson-CS} but includes
the constant mode (see (B.23) there) and the $\phi$ there should be
rescaled to $\phi/2\pi$ to agree with the definition here. So on putting
back the constant 
mode contribution in that work and changing the normalization of
$\phi$ we obtain
the (square root of the) determinant as
\be
\prod_{\alpha \in \mathbf{k}}  \left(\frac{\prod_{n}(2\pi n + i
  \alpha(\phi)/2\pi)}{i\alpha(\phi)/2\pi} \right)^{1-g} = \prod_{\alpha \in
\mathbf{k}_{+}}
\left(\frac{\sin^{2}{i\alpha(\phi)}/4\pi}{(i\alpha(\phi)/4\pi)^{2}}
\right)^{1-g} \equiv j_{\lg}(\phi)^{1-g}
\nonumber 
\ee
together with the phase (B.31)-(B.34) in \cite{Blau-Thompson-CS}. 
The integral over those $A^{\lt}$ which are non-constant leads to a
simple overall factor in front of the path integral. The integral over
the non-constant parts of the symplectic volume give rise to a
normalization which compensates that of the connections. The
calculation presumes $\phi$ constant (not just along the fibre) as in
the path integral over $\mathcal{A}$ it may be taken to be so.\qed

There are still 2 issues that we need to deal with:
\begin{enumerate}
\item Extend Proposition \ref{propdet} to general
  sections $\phi \in \Gamma(\Sigma, \mathrm{ad}\, P)$.
\item Make sure that supersymmetry is preserved.
\end{enumerate}

\begin{remark}{{\rm Both issues are resolved by recalling a basic
      tennet of renormalizable 
field theory: Upon regularising a theory it may be neccessary to add
to the Lagrangian local counterterms in order to restore symmetries
broken by the choice of regularization.}}
\end{remark}

We begin with the first issue. We know that $\Gamma(\phi)$ is
formally gauge invariant, 
\be
\Gamma(g^{-1}\,\phi\, g) = \Gamma(\phi), \;\;\;\; g \in \mathcal{G} \nonumber
\ee
However, the curvature 2-form $F_{A}$ that appears in the formula in
Proposition \ref{propdet} lies in the Cartan direction,
$F_{A}^{\lt}=dA^{\lt}$, and as it 
stands the result is only gauge invariant under gauge transformations
in the maximal torus, $g \in \mathrm{Map}(\Sigma,T)$. The source
for this is that the regularization adopted in \cite{Blau-Thompson-CS}
was only designed to preserve the Torus invariance.

It is straightforward to check that, in general, the absolute value of the
determinant is a function of $\phi^{2}$ and indeed agrees with the
function $\widehat{A}(\sqrt{\phi^{2}})$ so that this is invariant
under $\mathcal{G}$. Our difficulty, therefore,
rests with the phase and we perform a finite renormalization to put 
in the complete non-Abelian curvature 2-form which re-instates $\mathcal{G}$
invariance.

This is still not quite the end of the story. The one loop correction
is not supersymmetric. Or put another way we have not maintained the
original supersymmetry (\ref{susy}) at the level of the zero modes,
which is now,
\be
QA=i\psi, \;\;\; Q\psi=-d_{A}\phi, \;\;\; Q\phi=0 , \label{susy2}
\ee
We can add another finite renormalization 
\be
\exp{\frac{c_{\mathbf{g}}}{8\pi^{2}} \int_{\Sigma} \Tr{\psi \wedge
    \psi} }\nonumber
\ee
to correct this.

\begin{proposition}\label{1loopexact}{\rm{The gauge invariant and supersymmetric
      evaluation of the path integral along the fibres of $M$ is
\be
\exp{i\Gamma(A,\phi, \psi)} \simeq \exp{\left(i\frac{\pi}{2}\eta_{0}\right)}\,
\widehat{A}(i\phi) \wedge \exp{\left( \frac{c_{\mathbf{g}}}{4\pi^{2}}
    \int_{\Sigma} \Tr{(i \phi.F_{A}+ \frac{1}{2} \psi \wedge \psi +
     i \frac{p}{4\pi} \phi^{2}\omega)} 
  \right)} \nonumber
\ee 
}}
\end{proposition}
\begin{remark}{\rm{With
      the identifications that $\psi \simeq \Psi$ and $\phi
      \simeq -i\Phi$ and as $c_{\lg}=r$ we have
\be
\exp{i\Gamma(A,\phi, \psi)} \simeq\exp{\left(i\frac{\pi}{2}\eta_{0}\right)}
\,\mathrm{Todd}(\mathcal{A})\, \wedge \, \exp{\left( 
    -i \frac{p}{2\pi}c_{\lg} \Theta(\mathcal{A}) \right) } \nonumber
\ee
}}
\end{remark}

The path integral now becomes one over objects defined on the Riemann
surface directly and we have established Proposition \ref{onS}. \qed

\section{Wilson Lines and Parabolic Points}
How is the picture that we have obtained in the previous sections
affected by the inclusion of Wilson lines? Since our manifold is a
$S^{1}$ fibration there is a special 
class of knots which are located at point $x\in \Sigma$ on the base of
the fibration and which run along the fibre. To such a fibre knot we
associate
\be
W_{R}(x) = \Tr_{R}{P\exp{\left(\int_{S^{1}}\kappa \phi/2\pi\right)}} =
\Tr_{R}{\exp{\left(\phi(x)/2\pi\right)} } \nonumber
\ee
the second equality following from the condition that
$\iota_{\kappa}d\phi =0$. As the only addition to path integral
involves functions without 
dependence on the fibre, the calculation of the previous section goes
through unchanged.

A geometric way in which to add such traces is through the introduction
of co-adjoint orbits. Let $\lambda \in \lg^{*}$ ($\lg^{*}$ is the dual of $\lg$,
however, we identify the two so that an invariant inner product
$<f,\phi>\equiv \Tr{f\phi}$, $f \in \lg^{*}$, $\phi \in \lg$)
then the orbit through $\lambda$ is
\be
M_{\lambda}=\left\{ g^{-1}\lambda g;\; \forall g \in G \right\} \nonumber
\ee
while the stabilizer of $\lambda$ is
\be
G(\lambda) = \left\{ g \in G :\; g^{-1}\lambda g =\lambda\right\}\nonumber
\ee
If $\lambda$ is regular ($\det_{\mathbf{k}}{(\ad{\lambda})} \neq 0$) then
$G(\lambda)=T$ and we consider this case for now so that
$M_{\lambda}=G/G(\lambda)=G/T$. 

The homogeneous space $G/T$ comes equipped with a natural $G$
invariant symplectic 2-form (the Kirillov-Konstant form)
$\Omega_{\lambda}$ given by 
\be
\Omega_{\lambda}(X,Y) = <\lambda, \, [X,Y]> = \Tr{(\lambda\,[X,Y])
}\;\;\; X,\,Y\in \lg \nonumber
\ee
Kirillov tells us
that for $\lambda= \Lambda + \rho$ regular, $\Lambda$ an element of
the weight lattice and $\rho$ the Weyl vector then
\be
\Tr_{\lambda}{(\exp{\phi/2\pi})} = j_{\lg}^{-1/2}(\phi/2\pi)\,
\int_{M_{\lambda}} \exp{\left( i\frac{1}{2\pi}<\lambda, \phi> +
    \Omega_{\lambda}  \right)} \nonumber
\ee
Now we see that geometrically we should product in the co-adjoint
orbits so consider the space $\mathcal{A}_{3} \times \prod_{i} M_{R_{i}}$,
and we have
\be
Z_{k,G}[M,(x_{i},R_{i})] = \frac{1}{\mathrm{Vol}(\mathcal{G}_{3})}
\int_{\mathcal{A}_{3}\times \prod_{i} 
  M_{R_{i}}}  \exp{\left(I(\mathscr{A}) \right)}\, \prod_{i=1}
j_{\mathbf{g}}^{-1/2}(\phi(x_{i})/2\pi)\, \exp{\omega(M_{R_{i}})}\nonumber
\ee
where, in analogy with $\Omega(\mathcal{A})$,
\be
\omega(M_{R_{i}})= \frac{i}{2\pi}\Tr \lambda_{i}\phi(x_{i}) +
\Omega_{R_{i}}\nonumber 
\ee
We have the following:
\begin{lemma}{\rm{(Lemma 8.5 \cite{BGV}) The equivariant
      $\widehat{A}$-genus, $\widehat{A}_{\mathbf{g}}(X,G/T)$, of the
      Riemannian manifold $G/T$ and $j_{\mathbf{g}}^{-1/2}(X)$
      represent the same class in equivariant deRham cohomology.
}}
\end{lemma}
Consequently Proposition \ref{onSKnot} is proved. \qed

\begin{remark}{{\rm C. Beasley \cite{Beasley-Knot} has computed, in
      the spirit of \cite{Beasley-Witten}, the
      localization formula for
      $\left. Z_{k,G}[M,(x_{i},R_{i})]\right|_{\mathfrak{M}}$. This formula
    agrees with that in Proposition \ref{onSKnot} when restricted to
    $\mathfrak{M}$. 
}}
\end{remark}

\bibliographystyle{amsplain}


\end{document}